\documentclass[11pt]{amsart}
\usepackage[cp1251]{inputenc}

\usepackage{amssymb,amsmath,amscd,color}

\theoremstyle{plain}
\newtheorem{thm}{Theorem}
\newtheorem{theorem}[thm]{Theorem}

\newtheorem{prop}{Proposition}

\theoremstyle{definition}

\theoremstyle{remark}

\newtheorem*{remark*}{Remark}

\newcommand{\cL}{{\mathcal{L}}}

        \newcommand{\field}[1]{{\mathbb{#1}}}
        \newcommand{\NN}{\field{N}}
        \newcommand{\ZZ}{\field{Z}}
        
        \newcommand{\RR}{\field{R}}
        
        \newcommand{\HH}{\field{H}}

\allowdisplaybreaks

\begin{document}

\title[Trace formula for the magnetic Laplacian]{Trace formula for the magnetic Laplacian on a compact hyperbolic surface}

\author[Y. A. Kordyukov]{Yuri A. Kordyukov}\address{Institute of Mathematics, Ufa Federal Research Centre, Russian Academy of Sciences, 112~Chernyshevsky str., 450008 Ufa, Russia} \email{yurikor@matem.anrb.ru}

\author[I. A. Taimanov]{Iskander A. Taimanov}\address{Sobolev Institute of Mathematics, 4 Acad. Koptyug avenue, and Novosibirsk State University, Pirogova st 1, 630090, Novosibirsk, Russia}\email{taimanov@math.nsc.ru}

\thanks{The second author (I.A.T.) was partially supported by  the Mathematical Center in Akademgorodok under the agreement No. 075-15-2022-282 with the Ministry of Science and Higher Education of the Russian Federation.}



\begin{abstract}
We compute the trace formula for the magnetic Laplacian on a compact hyperbolic surface of constant curvature with constant magnetic field for energies above the Mane critical level of the corresponding magnetic geodesic flow. We discuss the asymptotic behavior of the coefficients of the trace formula 
when the energy approaches the Mane critical level.  
\end{abstract}

\dedicatory{To the memory of Alexey Borisov}

\date{}

 \maketitle
\section{Introduction}

We consider a classical mechanical system which describes the motion of a charged particle in an external magnetic field on a Riemannian manifold. In addition to the dynamical and variational problems for such flows, which have been intensively studied over the last years, there are many interesting questions concerning the relations between the classical dynamics and the spectral properties of the corresponding quantum Hamiltonian, which is given by the magnetic Laplacian.  Recently we have addressed a few of them in \cite{KT1,KT2}. In \cite{KT2} we constructed the quasi-classical approximation for the eigenfunctions of the magnetic Laplacians. This article continues the research started in \cite{KT1}.

We study the Guillemin--Uribe trace formula for magnetic geodesic flows on hyperbolic surfaces with a constant magnetic field and on sufficiently low energy levels. The dynamics in this situation is mostly determined by the ratio $\sqrt{E_0}/B$ where $E_0 = |p|^2$ is the squared norm of the momentum and the magnetic field takes the form $B\, d\mathrm{vol}$ where $d\mathrm{vol}$ is the area form corresponding to the hyperbolic metric. Without loss of generality, we assume that $B=1$. Then the level
$$
E_0=\mu_0 : = 1
$$
is known as the Mane critical level \cite{CMP,CFP} for this system which is integrable for $E_0<1$ and chaotic for $E_0>1$.

The magnetic geodesic flow on a Riemannian manifold $(M,g)$ is determined by a magnetic field $F$, which is a closed $2$-form. The magnetic Laplacian is defined iff $F$ satisfies the integrality condition $\left[\frac{F}{2\pi}\right] \in H^2(M;\ZZ)$. In this case one can define the Hermitian line bundle $L$ on $M$ with a Hermitian connection such that $F$ is the curvature of this connection and the family of the magnetic Laplacians $\Delta^{L^N}, N\in \NN,$ acting on sections of its tensor powers $L^N$.

Let us denote by $\nu_{N,j}, j=0,1,2,\ldots$, the eigenvalues of $\Delta^{L^N}$ taken with multiplicities
and put $\lambda_{N,j}=\sqrt{\nu_{N,j}+N^2}$.
Given a constant $E>1$ and an arbitrary function $\varphi\in \mathcal S(\RR)$, let us introduce the sequence
\[
Y_N(\varphi)=\sum_{j=0}^{\infty}\varphi(\lambda_{N,j}-EN),\quad N\in \mathbb N.
\]
The Guillemin--Uribe trace formula  \cite{Gu-Uribe89} describes the asymptotic expansion of $Y_N$ as $N \to \infty$ in terms of the magnetic geodesic flow on the energy surface $E_0 = |p|^2 = E^2-1$ under some additional assumption (the cleanness of the flow):
\[
Y_N(\varphi)\sim \sum_{k=0}^\infty c_k(N,\varphi)N^{\dim M -1-k},\quad N\to \infty,
\]
where $c_k(N,\varphi)$ are bounded in $N$. Considered as functionals of the Fourier transform $\hat \varphi$, the coefficients $c_k$ are distributions supported in the union of zero and the period set of closed magnetic geodesics. They are rapidly oscillating in $N$ and the frequencies of these oscillations are given by the actions of closed magnetic geodesics. The contribution of zero is often called the Weyl term, because it is related with the asymptotic formula for the eigenvalue distribution function.

In \cite{KT1} we considered the constant magnetic fields on compact hyperbolic surfaces and computed the first two coefficients $c_0$ and $c_1$ of this expansion for the energy levels below the Mane critical level, i.e., for $E_0 <1$ (Theorem 2). Here we do that for $E_0\geq 1$ (Theorem 3). 

The coefficients of the trace formula depend on the energy level $E_0$ as follows:

\begin{itemize}
\item
for $E_0<1$, i.e., below the Mane critical level, the classical dynamics is integrable, the flow is periodic, and periodic trajectories are lifted to hyperbolic circles on the universal covering. The periodic orbits form two-dimensional Liouville tori which contribute to the leading coefficient $c_0$;

\item
for $E_0=1$, i.e., on the Mane critical level, the classical dynamics is given by the horocyclic flow which has no nontrivial periodic orbits. Therefore, the trace formula reduces to the Weyl term. We can also observe the convergence of the contribution of closed magnetic geodesics to $0$ as $E_0\to 1$. This happens because the periods of primitive closed magnetic geodesics go to infinity as $E_0\to 1$ (both from below and from above);

\item 
for $E_0>1$, i.e., above the Mane critical level, the magnetic geodesic flow is chaotic. The closed magnetic geodesics are noncontractible, nondegenerate, and isolated. They don't contribute to the leading coefficient $c_0$, which coincides with the Weyl term in this case. The magnetic system looks similar to the system, which describes the motion of a free particle on the surface. This can be seen both at the classical and quantum level. At the classical level, the magnetic geodesic flow for the energy level $E_0$ is conjugated to the Riemannian geodesic flow. At the quantum level, there is a relation between the spectrum  of the magnetic Laplacian $\Delta^{L^N}$ on the half-line $(N^2,\infty)$ and the spectrum of the Laplace-Beltrami operator. One should note that the coefficient at $N^2$ here is exactly the Mane critical level $\mu_0=1$. We use these facts to give an alternative computation of the trace formula, reducing it to the case of the Laplace-Beltrami operator.
\end{itemize}
 
In \cite{Taimanov2004} the second author (I.A.T.) showed how to construct an additional real-analytic first integral for this system below the Mane critical level. Regretfully due to the brevity of the original communication the scenario of the degeneration of these integrals as the energy tends to the Mane cirtical level was skipped. For completeness we expose it here in Theorem 1.

\section{Classical system}

Let $M^2$ be a surface with metric of constant negative curvature $K$. This means that it is isometric to 
the quotient of the space $\HH$ with respect to some discrete group $\Gamma$ which acts by isometries. 
 
We consider two models of $\HH$ which are helpful for different reasons:

1) $\HH$ is the upper half-plane $\{(x,y) \in \RR^2\,:\, y>0\}$
endowed with the metric
\begin{equation}
\label{metric}
g = \frac{dx^2 + dy^2}{y^2}.
\end{equation}
The curvature of this space is equal to
$$
K = - 1
$$
and the full isometry group is $PSL(2,\RR) = SL(2,\RR)/\pm 1$ which acts by fractional linear transformations:
$$
z = x+iy \to \frac{az+b}{cz+d}, \ \ \ \det \left(\begin{array}{cc} a& b \\ c & d \end{array}\right) = 1.
$$

2) Let us take the $(1+2)$-dimensional Minkowski space $\RR^3_{1,2}$ endowed with the metric
$$
ds^2 = dt^2 - dx^2 - dy^2
$$
and consider in it the pseudosphere
$$
t^2 - x^2 - y^2 = 1
$$ 
endowed with the induced metric (multiplied by $-1$). In the pseudospherical coordinates
$\rho,r,\varphi$ such that
$$
t = \rho \cosh r, \ \ x = \rho \sinh r \cos \varphi, \ \ y = \rho \sinh r \sin \varphi
$$
the pseudosphere and the induced metric are as follows:
$$ 
\rho = 1, \ \ ds^2 = dr^2 + \sinh^2 r d\varphi^2.
$$ 
The full isometry group of the pseudosphere is $O_0(1,2)$, the connected component of the unity in the group $O(1,2)$ of pseudoorthogonal linear transformations of $\RR^3_{1,2}$.
Therefore, $r$ and $\varphi$ are coordinates on the pseudosphere and the mapping
\begin{equation}
\label{map}
(r,\varphi) \to z = i\frac{1 - w}{1+ w}, \ \ \ w  = \frac{\sinh r}{1+\cosh r}e^{i\varphi},
\end{equation}
establishes an isometry of the pseudosphere to the upper half-plane with the metric (\ref{metric}).
 
By a constant magnetic field on a surface we mean a two-form $F = B d\mathrm{vol}$ which is a constant multiple of the volume form. In our case
\begin{equation}
\label{def:F}
F = B \frac{dx \wedge dy}{y^2} \ \ \ \mbox{or} \ \ \ F =  B \sinh r dr \wedge d\varphi,
\end{equation}
where $B$ is a constant. 

\subsection{The Lagrangian formalism}

The motion of a charged particle in an arbitrary magnetic field is described by the Euler--Lagrange equations for the lagrangian
$$
\cL(x,\dot{x}) = \frac{|\dot{x}|^2}{2} + A_i\dot{x}^i, \ \ \ x \in M, \, \dot{x} \in T_xM,
$$
where $A= A_idx^i$ is the vector potential of the magnetic field. In our case we take
\begin{equation}
\label{def:A}
A = \frac{B}{y}dx \ \ \ \mbox{or} \ \ \ A = B \cosh r d\varphi.
\end{equation}

We prefer to use the pseudosphere model for calculations. 
The lagrangian takes the form
$$
\cL = \frac{1}{2}(\dot{r}^2 + \sinh^2 r \dot{\varphi}^2) + B \cosh r \dot{\varphi}.
$$
The Euler--Lagrange equations for this lagrangian are as follows:
$$
\ddot{r} = \sinh r\, \dot{\varphi} (\cosh r\, \dot{\varphi} + B), \ \ \frac{d}{dt}(\sinh^2 r \,\dot{\varphi} + B \cosh r)=0.
$$
We have two first integrals: the kinetic energy
$$
\frac{E_0}{2} = \frac{1}{2}(\dot{r}^2 + \sinh^2 r \dot{\varphi}^2),
$$
and, since $\frac{\partial \cL}{\partial \varphi} = 0$, the  momentum with respect to $\varphi$
\[
I=\sinh^2 r \,\dot{\varphi} + B \cosh r.
\]

Since $O_0(1,2)$ acts on the pseudosphere transitively, given a trajectory $c(t)=(r(t),\varphi(t))$,
we may assume that
$$
\dot{r}=0 \ \ \mbox{for $t=0$}.
$$  
If in addition
$$
\cosh r \, \dot{\varphi} + B = 0,
$$
then $\ddot{r}=0$ and the trajectory has the form
$$
r = \mathrm{const}, \ \ \varphi = -\frac{B}{\cosh r}t + \varphi_0, \ \varphi_0 = \mathrm{const}.
$$
It is easy to calculate that in this case
$$
\frac{E_0}{2}  = \frac{B^2 \tanh^2 r}{2}.
$$ 
Therefore we derive

\begin{prop}
If $0 < E_0 <  B^2$, then all trajectories are circles of radius
$$
R_{E_0} = \frac{1}{2} \log \frac{1+u}{1-u} \ \ \ \mbox{with $u = \frac{\sqrt{E_0}}{B}$},
$$
i.e.
$$
\tanh R_{E_0} = \frac{\sqrt{E_0}}{B}.
$$
\end{prop}

{\sc Remark.} For surfaces of constant positive curvature $K = 1$ analogous calculations show
that for every energy  $E_0$ all trajectories are circles of radius
$$
R_{E_0} =  \arctan \frac{\sqrt{E_0}}{B}.
$$

{\sc Hyperbolic cirles ($E_0 < B^2$).} By (\ref{map}), it is shown that the circles of radius
$R$ centered at the origin are mapped into Euclidean circles (on the upper-half plane with the metric (\ref{metric})) of radius 
$r = \sinh R$ centered at $z = i\cosh R$.
Since the hyperbolic circles are invariant under the actions of
$$
\left( \begin{array}{cc} \lambda^{1/2} & a \lambda^{-1/2} \\ 0 & \lambda^{-1/2}\end{array}\right) \in SL(2,\RR),
$$
where $\lambda$ is real and positive and $a \in \RR$, we conclude that 
all hyperbolic circles of radius $r$ are realized
by Euclidean circles with center at $z$ and radius $\rho$ where
\begin{equation}\label{e:circle}
z = i \lambda \cosh R +a, \ \ \  \rho = \lambda \sinh R, \ \ \lambda \in \RR^+, a \in \RR.
\end{equation}
By simple computation, it can be shown that the geodesic curvature $\varkappa$ of the hyperbolic circle of radius $R$ meets the equality
$$
\varkappa^2  = \frac{1}{\tanh^2 R} = \frac{B^2}{E_0} > 1.
$$

{\sc Horocycles ($E_0 = B^2$).} The formula \eqref{e:circle} in the limit 
$$
r \to \infty, \ \ \lambda \to \frac{\mu}{\sinh R}
$$
gives us the horocycles which are realized 
by Euclidean circles with center at $z_0$ and radius $\rho_0$ where
$$
z_0 = i \mu  +a, \ \ \  \rho_0 = \mu, \ \ \mu \in \RR^+, a \in \RR.
$$
To complete the description of horocycles, we have to add the images of these circles by the action of $PSL(2,\RR)$ which are Euclidean lines
$$
y = \mathrm{const} > 0.
$$
They correspond to the case when the horocycle touches $\{y=0\}$ at the infinite point
$a = \infty$. For the geodesic curvature $\varkappa$ of a horocycle we have
$$
\varkappa^2 = \frac{B^2}{E_0} = 1.
$$

\subsection{The Hamiltonian formalism}

In the Hamiltonian formalism, the motion of a charged particle in the magnetic field \eqref{def:F} is described by the magnetic geodesic flow, which is the  Hamiltonian flow on the phase space $X=T^*M$ equipped with the twisted symplectic form 
\[
\Omega=dp_x\wedge dx+dp_y\wedge dy+\frac{B}{y^2}dx\wedge dy.
\]
We have to remark that we consider two different Hamiltonian functions:
the kinetic energy
$$
H_0 = \frac{y^2(p_x^2+p_y^2)}{2} = \frac{|p|^2}{2}
$$ 
and the Hamiltonian  
$$
H = (2H_0+1)^{1/2}=\left(y^2(p_x^2+p_y^2)+1\right)^{1/2},
$$
which is more relevant for our considerations of the trace formula, and denote by $E_0/2$ and $E$ the values of $H_0$ and $H$, respectively.
There is a relation
$$
E_0 = E^2 -1.
$$
The corresponding Hamiltonian flows are related by time reparameterization.

The Hamiltonian system defined by $H$ has the form:
\begin{equation}\label{e:1.6}
\dot x=\frac{y^2}{H}p_x,
\quad \dot y=\frac{y^2}{H}p_y,
\quad \dot p_x=\frac{B}{H}p_y,
\quad \dot p_y=-\frac{y}{H}(p_x^2+p_y^2)-\frac{B}{H}p_x.
\end{equation}

Denote 
\begin{equation}\label{e:XE}
X_E=H^{-1}(E)=\{y^2(p_x^2+p_y^2)=E^2-1\}. 
\end{equation}
It is easy to see that $E>1$ is a regular value of $H$, and, therefore, $X_E$ is a smooth submanifold of  $T^*M$. The reduction of the system \eqref{e:1.6} to $X_E$ is given by 
\[
\dot x=\frac{y^2}{E}p_x,
\quad \dot y=\frac{y^2}{E}p_y,
\quad \dot p_x=\frac{B}{E}p_y,
\quad \dot p_y=-\frac{y}{E}(p_x^2+p_y^2)-\frac{B}{E}p_x.
\]
Let us introduce on 
$X_E= \{y^2(p_x^2 + p_y^2) =  E_0 \}$
the coordinates
$(x,y,\theta)$: 
\[
p_x= \frac{\sqrt{E_0}}{y}\cos\theta, \quad p_y= \frac{\sqrt{E_0}}{y}\sin\theta
\]
in which the system takes the form
\begin{equation}
\label{hom}
\dot x=\frac{\sqrt{E_0}}{E} y\cos\theta,
\quad \dot y=\frac{\sqrt{E_0}}{E} y\sin\theta,
\quad \dot \theta=-\frac{\sqrt{E_0}}{E}\cos\theta - \frac{B}{E}.
\end{equation}

This system has an evident conservation law:
$$
f = \frac{\dot{\theta}}{y}, \ \ \ \frac{df}{dt} = 0.
$$
Therefore
$$
\dot{\theta} = cy \ \ \mbox{for some constant $c$}.
$$
If $c=0$ then the system has a solution whose trajectory is a Euclidean line.
This is possible if
$$
\frac{\sqrt{E_0}}{E}\cos\theta + \frac{B}{E} = 0
$$ for some $\theta$.
We rewrite the last equality as
\begin{equation}
\label{angle}
\cos \theta = - \frac{B}{\sqrt{E_0}}.
\end{equation}
Since $-1 \leq \cos \theta \leq 1$, such a $\theta$ exists if and only if
$$
E_0 \geq B^2.
$$
The case $E_0 = B^2$ corresponds to horocycles and the remaining case to  hypercycles.

{\sc Hypercycles ($E_0 > B^2$)} are given by the Euclidean lines which meet the axis $\{y=0\}$ at the constant angle $\theta$ satisfying (\ref{angle}) and the images of these lines under the action of $PSL(2,\mathbb R)$.
Let us take such a line and consider another line which meets the axis $\{y=0\}$ at the same point and is orthogonal to it. It is a trajectory of the geodesic flow. It is easy to notice that the hypercycle and this geodesic are equidistant. Therefore for any hypercycle there exists an equidistant geodesic, the distance is the same for all hypercycles, and these two flows are conjugate after a constant time scaling.  
The formula for the geodesic curvature $\varkappa$ of hypercycles looks the same as for hyperbolic cycles and horocycles:
$$
\varkappa^2 = \frac{B^2}{E_0} <1.
$$ 

\subsection{The Lie group approach}\label{s:Lie}

Let us consider the ANK decomposition of the group  $PSL(2,\mathbb R)$:
\[
\begin{pmatrix}
y^{1/2} & 0\\
0 & y^{-1/2}
\end{pmatrix}
\begin{pmatrix}
1 & x\\
0 & 1
\end{pmatrix}
\left[
\begin{pmatrix}
\cos \frac{\varphi}{2} & \sin \frac{\varphi}{2}\\
-\sin \frac{\varphi}{2} & \cos \frac{\varphi}{2}
\end{pmatrix}\big/
\begin{pmatrix}
\pm 1 & 0 \\
0 & \pm 1
\end{pmatrix}
\right],
\] 
where $x, y \in \mathbb R$ and $y >0$.
It gives a unique representation of an element from $PSL(2,\mathbb R)$ as a product of elements from
the canonical subgroups $A$, $N$, and $K$. Since the inversion $g \to g^{-1}$ maps the ANK decomposition into the KNA decomposition and the products of subgroups $AN$ and $NA$ coincide, we 
have the canonical KAN decomposition which is also known as the Iwasawa decomposition. 

The KNA decomposition was used in \cite{GF} for describing the geodesic flow on $\mathbb H$.
There is the isomorphism 
\[
PSL(2,\mathbb R) \to S\mathbb H\cong  \{(x,y,p_x,p_y) : y^2(p_x^2+p_y^2)= 1\} 
\]
of the form
\[
(x,y,\varphi) \to \left( x,y, p_x = \frac{1}{y} \cos \left(\varphi +\frac{\pi}{2}\right), 
p_y = \frac{1}{y} \sin \left(\varphi +\frac{\pi}{2}\right)\right) , 
\] 
i.e., $\theta = \varphi+\frac{\pi}{2}$.

The geodesic which starts at $i \in \mathbb H$ and is directed along the imaginary axis is 
the orbit of $i$ under the action of $A$. Hence, all geodesic orbits are uniquely parametrized 
by elements from $KN \subset PSL(2,\mathbb R)$. 

The Lie algebras of the one-dimensional subgroups $A, N$, and $K$ are generated by
\[
e_1 = 
\begin{pmatrix}
1/2 & 0 \\
0 & -1/2
\end{pmatrix}, \ \ 
e_2= 
\begin{pmatrix}
0 & 1\\
0 & 0
\end{pmatrix}, \ \
e_3 = \begin{pmatrix}
0 & 1/2 \\
-1/2 & 0
\end{pmatrix}.
\]
It is clear that under the isomorphism  $S\mathbb H\cong  PSL(2,\mathbb R)$
the left-invariant vector field corresponding to $e_3$ is
\[
V_3 = 
\frac{\partial}{\partial \varphi}
\]
and such a field corresponding to $e_1$ is
\[
V_1 = - y \sin \varphi \frac{\partial}{\partial x} + y \cos \varphi \frac{\partial}{\partial y} + 
\sin \varphi \frac{\partial}{\partial \varphi}.
\]

Let us rewrite (\ref{hom}) as the equations on $S\mathbb H$:
\begin{equation}\label{SH}
 \dot x=- \alpha y\sin\varphi,\quad \dot y=\alpha y\cos\varphi,\quad \dot \varphi=\alpha \sin \varphi +\beta
\end{equation}
with 
\begin{equation}\label{ab}
\alpha=\frac{\sqrt{E_0}}{E}, \quad \beta=- \frac{B}{E}, 
\end{equation}
and notice that they describe the motion along trajectories of the left-invariant vector field 
\[
\alpha V_1 + \beta V_3.
\]

Denote by $\Phi^{\alpha,\beta}_t$ the flow given by the right translation by $\exp(t(\alpha E_1+\beta E_3))$. 
Since
\[
\det (\alpha V_1 + \beta V_3)=-1/4(\alpha^2-\beta^2) = -\frac{1}{4}\frac{E_0 - B^2}{E^2},
\]
we have to specialize three cases when $\det >0, \det = 0$, and $\det <0$:

\begin{enumerate}
\item $\det >0$, i.e., $E_0 < B^2$: the hyperbolic cycles;

\item $\det =0$, i.e., $E_0 = B^2$: the horocycle flow;

\item $\det <0$, i.e., $E_0>B^2$: the hypercycle flow. In this case  the flow $\Phi^{\alpha,\beta}_t$ is conjugate to the flow $\Phi^{\delta,0}_t$ for certain $\delta$. For $B=1$ the condition $\det  < 0$ is equivalent to
\[
E>\sqrt{2} \ \ \ \mbox{or} \ \ \ E_0 >1.
\]
\end{enumerate}

The constant $\delta$ can be found from the condition
\[
\det (\alpha V_1 + \beta V_3)=\det (\delta V_1)=-\frac{1}{4}\delta^2,
\]
which gives $\delta=\frac{\sqrt{E_0 - B^2}}{E}$. This approach to describing magnetic geodesic flows was initiated in \cite{Sunada} (see also \cite{Butler,CFP}).

\subsection{The (Mane) critical level}\label{s:Mane}

The qualitative behavior of the magnetic geodesic flow on $M$ depends on the ratio $\tau = \frac{B^2}{E_0}$. Since the flow on $\mathbb H$ is completely integrable, we see from the exact formulas for trajectories that the dynamics is different for $\tau <1, \tau =1$ and $\tau>1$. 

At the beginning of 1960s just after the emergence of the Kolmogorov entropy theory Arnold had shown that, given $B=1$, if the metric entropy (with respect to the Liouville measure) of the geodesic flow is equal to $h(0)$, then the metric entropy of the hypercycle flow is equal to $h(\varkappa) = h(0)\sqrt{1-\varkappa^2}$ 
and it vanishes for $\varkappa^2 \geq 1$ \cite{Arnold}.

The relation of this flow of linear elements to magnetic fields was not discussed in \cite{Arnold}
and probably it was first considered by Ginzburg \cite{Ginzburg} who pointed out that the horocycle flow 
on a closed hyperbolic surface gives an example of a magnetic geodesic flow without periodic trajectories.
Until recently this is the only known such an example.  

A systematic study of magnetic geodesic flows started in the early 1980s \cite{Novikov1}. Therewith, such flows appeared as reductions of mechanical systems (the Kirchhoff equation, mechanical tops) and their periodic trajectories do not describe motions of charge particles in real magnetic fields. Although great progress was achieved in the study of the periodic variational problem for such systems, the original periodic problems for explicit mechanical systems mostly remain unsolved \cite{Novikov2}.

In \cite{Taimanov2004} the second author (I.A.T.) mentioned that for $\tau>1$ the flow is easily integrable and its additional (to the kinetic energy) first integral $F$ can be constructed from any smooth function $f: M \to {\mathbb R}$ on the hyperbolic surface $M$. Indeed, for every point $q=(p.\xi) \in SM$ we consider the hyperbolic circle $\gamma$ on $M$ such that $\gamma(0)=p,\dot{\gamma}(0)=\xi$, take the center $c_\gamma$ of this circle and put   
$$
F(q) = f(c_\gamma).
$$
These integrals were successfully used in \cite{BNS} for a quantization of periodic magnetic geodesics.

As $\tau \to 1$ the integrability disappears and for $\tau<1$ the flow is chaotic. 

Due to the brevity of a short communication  \cite{Taimanov2004}, the scenario of the degeneration of the first integral $F$ was skipped and we describe it here.

Let us assume that $E_0=1$. 
Take 
\[
p \in M, \xi \in T_p M, |\xi|=1, q=(p,\xi),
\]
\[
\eta \in T_p M  \ \mbox{such that $\xi \perp \eta, \xi\wedge \eta <0$}.
\]
Let us draw the geodesic 
$\gamma_{p,\xi}: [0,\infty)\to M$ with the inital data
$\gamma_{p,\xi}(0) = p, \dot{\gamma}_{p,\xi}(0) = \xi$.

\begin{theorem}
Given $E_0 =1 $, $B>1$, and a smooth function $f: M \to {\mathbb R}$, we have the first integral
$F_B$:
$$
F_B(p,\xi) = f\left(\gamma_{p,\xi}\left(\frac{1}{2}\log \left(\frac{B+1}{B-1}\right)\right)\right).
$$
Regular contours lying on an energy level $F_B = \mathrm{const}$ give rise to invariant two-dimensional tori. As $B \to 1$ the first integral $F_B$ degenerates.
\end{theorem}

From the modern point of view this example is considered as a particular case of the Mane critical level 
\cite{CMP,CFP}. Given $B=1$, the energy level $E_0=1$ is (Mane) critical.  We skip its definition here
however in the rest of the article we study how  transition through the critical level affects the (Guillemin--Uribe) trace formula for the magnetic geodesic flow. 

\section{Quantum system and trace formula}

\subsection{Quantum Hamiltonian}\label{s:spectrum}
To quantize a classical magnetic system on a Riemannian manifold $(M,g)$ defined by a magnetic field $F$, it is necessary that the de Rham cohomology class of the form $\frac{1}{2\pi}F$ is integral:
 $$
 \left[\frac{1}{2\pi} F\right] \in H^2(M;\ZZ).
 $$
In this case, it is the first Chern class of a line bundle $L$ on $M$:
 $$
 c_1(L) = \left[\frac{1}{2\pi} F\right],
 $$
and the quantum Hamiltonian $\Delta^L$ (the magnetic Laplacian) acts on sections of $L$. Its definition depends on the choice of a Hermitian connection $\nabla^L$ on $L$ (a vector potential of the magnetic field). 

For the upper half-plane model of the hyperbolic plane $\mathbb H$ endowed with metric \eqref{metric} and constant magnetic field \eqref{def:F} with an arbitrary $B$, we can take the Hermitian line bundle $\tilde L$ to be trivial and the connection form of the connection $\nabla^{\tilde L}$ on $\tilde L$ to be given by \eqref{def:A}. The corresponding magnetic Laplacian on $\mathbb H$ is given by
\[
\Delta^B=-y^2\left(\left(\frac{\partial}{\partial x}-iBy^{-1}\right)^2+\frac{\partial^2}{\partial y^2}\right).
\] 
Such an operator first appeared in the theory of automorphic forms, where it is known as the Maass Laplacian. 
More precisely, it is related with the Maass Laplacian
\[
D_B=y^2\left(\frac{\partial}{\partial x^2}+\frac{\partial^2}{\partial y^2}\right)-2iBy\frac{\partial}{\partial x}
\]
by the formula 
\[
\Delta^B=-D_B+B^2.
\]
The relation between magnetic trajectories on the hyperbolic plane and the corresponding spectral properties of the magnetic Laplacian was first discussed in \cite{Comtet87} where it was also observed that the magnetic Laplacian on the hyperbolic plane is given by the Maass operator.

For a compact hyperbolic surface $M=\Gamma\setminus\mathbb H$, where $\Gamma\subset PSL(2,\mathbb R)$ is some discrete group of isometries, endowed with metric \eqref{metric} and constant magnetic field \eqref{def:F}, the Hermitian line bundle $L$ exists iff the quantization condition 
\[
(2g-2)B\in \mathbb Z
\]
holds true, where $g$ is the genus of $M$.

Under this assumption, we can choose a Hermitian line bundle $L^B$ on $M$ such that its smooth sections on $M$ are identified with smooth functions $\psi$ on $\mathbb H$, satisfying the condition
\begin{equation}\label{e:sections}
\psi(h z)=\psi(z)\exp(i2B\arg(cz+d))=\left(\frac{cz+d}{|cz+d|}\right)^{2B}\psi(z)
\end{equation}
for any $z\in \mathbb H$ and $h=\begin{pmatrix} a & b\\ c & d\end{pmatrix}\in \Gamma$.  

We will fix $B=1$ and denote by $L$ the Hermitian line bundle $L^B$ for $B=1$. For any $N\in \mathbb N$, the bundle $L^N$ is the $N$th tensor power of $L$, $L^N=L^{\otimes N}$, and the space $C^\infty(M,L^N)$ of its smooth sections is identified with the space $\mathfrak F_N$ of smooth functions $\psi$ on $\mathbb H$, satisfying the condition \eqref{e:sections} with $B=N$. The parameter $N$ plays the role of a semiclassical parameter for the symplectic manifold $(M,F)$, and the limit $N\to \infty$ can be considered as the semiclassical limit. One can show that $D_N : \mathfrak F_N\to \mathfrak F_N$. Therefore, the magnetic Laplacian $\Delta^{L^N}$ on $C^\infty(M,L^N)$ corresponds to the restriction of the operator $-D_N+N^2$ to $\mathfrak F_N$ under isomorphism $C^\infty(M,L^N)\cong \mathfrak F_N$. 

The spectrum of $\Delta^{L^N}$ is computed by means of the Maass operators \cite{Maass53}, which are first order differential operators on $\mathbb H$ given by
\[
K_N=(z-\bar z)\frac{\partial }{\partial z}+N=2iy^{1-N}\frac{\partial }{\partial z}y^N,
\]
\[
L_N=(\bar z- z)\frac{\partial }{\partial \bar z}+N=-2iy^{1+N}\frac{\partial }{\partial \bar z}y^{-N}.
\]
Recall some basic properties of these operators: 
\[
\overline{K_N}=L_{-N}, \quad K_N^*=-L_{N+1}.  
\]
\[
D_N=L_{N+1}K_N+N(N+1)=K_{N-1}L_N+N(N-1). 
\]
\[
D_{N+1}K_N=K_ND_N, \quad D_{N}L_{N+1}=L_{N+1}D_{N+1}. 
\]
\[
K_N : \mathfrak F_N\to \mathfrak F_{N+1}, \quad L_N : \mathfrak F_N\to \mathfrak F_{N-1}.
\]
For the magnetic Laplacian $\Delta^{L^N}$, we have
\[
\Delta^{L^N}=K_N^*K_N-N=L_{N}^*L_N+N.
\]

Using the Maass operators, one can compute the spectrum of $\Delta^{L^N}$ on the interval $[0,N^2]$ \cite{Roelcke} (see also \cite{Elstrodt,Comtet-Houston85, FV}). It consists of eigenvalues 
\begin{equation}\label{e:nu-i}
\nu^{(i)}_{N,k}= (2k+1)N-k(k+1), \quad 0\leq k\leq N-1,
\end{equation}
with multiplicity
\[
m_{N,k}=(g-1)(2N-2k-1), \quad 0\leq k\leq N-1.
\] 

The theory of Maass operators also allows us to relate the spectrum of $\Delta^{L^N}$ on the half-line $(N^2, \infty)$ with the spectrum of the Laplace-Beltrami operator on $M$. More precisely, let $\Delta_{\mathbb H}$ denote the Laplace-Beltrami operator of the metric \eqref{metric} on $\mathbb H$:
\[
\Delta_{\mathbb H}=-y^2\left(\frac{\partial}{\partial x^2}+\frac{\partial^2}{\partial y^2}\right)=-D_0
\]
and $\Delta_M$ denote the Laplace-Beltrami operator on $M$. The operator $\Delta_M$ on $C^\infty(M)$ corresponds to the restriction of the operator $\Delta_{\mathbb H}=-D_0$ to $\mathfrak F_0$ under isomorphism $C^\infty(M)\cong \mathfrak F_0$. Denote by 
\[
\lambda_0=0<\lambda_1\leq \lambda_2\leq \ldots, \lambda_\ell\to +\infty,
\] 
the eigenvalues of $\Delta_M$ (taking into account multiplicities):
\[
\Delta_M \psi_\ell=\lambda_\ell \psi_\ell, \quad \psi_\ell \in C^\infty(M)\cong \mathfrak F_0.
\]
 
Using the properties of the Maass operators, one can show (see, for instance, \cite[p. 146]{Fay77}) that the functions
\[
u_{N,\ell}=\frac{1}{c_{N,\ell}}K_{N-1}\ldots K_0\psi_\ell\in \mathfrak F_N,
\]
with some suitable constants $c_{N,\ell}$, are orthonormal eigenfunctions of $D_N$ with the eigenvalue $-\lambda_\ell$. It follows that the eigenvalues of the magnetic Laplacian $\Delta^{L^N}$ on the half-line $(N^2, \infty)$ are given by
\begin{equation}\label{e:nu-c}
\nu^{(c)}_{N,\ell} = \lambda_\ell+N^2, \quad \ell=0,1,2,\ldots.
\end{equation}

\subsection{The trace formula}
Let $(M,g)$ be a compact Riemannian manifold equipped with a magnetic field $F$, satisfying the integrality condition, $L$ the associated Hermitian line bundle on $M$ with Hermitian connection and $\Delta^{L^N}, N\in \NN,$ the magnetic Laplacian, acting on sections of $L^N$. Denote by $\{\nu_{N,j}, j=0,1,2,\ldots\}$ the eigenvalues of $\Delta^{L^N}$ taken with multiplicities. Put
\begin{equation}\label{e:def-lambda}
\lambda_{N,j}=\sqrt{\nu_{N,j}+N^2}.
\end{equation}
Fix $E>1$. For an arbitrary function $\varphi\in \mathcal S(\RR)$, we introduce the sequence
\begin{equation}\label{e:Yp}
Y_N(\varphi)=\sum_{j=0}^{\infty}\varphi(\lambda_{N,j}-EN),\quad N\in \mathbb N.
\end{equation}
The Guillemin-Uribe trace formula \cite{Gu-Uribe89} describes the asymptotic expansion, as $N\to\infty$, of the sequence $Y_N$ given by \eqref{e:Yp} with some $E>1$ and $\varphi\in \mathcal S(\RR)$ with compactly supported Fourier transform in terms of the magnetic geodesic flow on the energy level $X_E$ (see \eqref{e:XE}) under the assumption on the flow to be clean. A survey of basic notions and results related with the Guillemin-Uribe trace formula is given in \cite{KT1}. In \cite{KT1}, we have also provided some concrete examples of its computation. In particular, we computed the trace formula in the current setting of hyperbolic surfaces with constant magnetic fields in the case $1<E<\sqrt{2}$. We note that the threshold value $E=\sqrt{2}$ corresponds exactly to the Mane critical level $\mu_0=1$ discussed above. Let us recall the result. 

Let $M=\Gamma\setminus\mathbb H$ be the compact hyperbolic surface endowed with metric \eqref{metric} and constant magnetic field \eqref{def:F} with $B=1$ and let $L=L^1$ be the Hermitian line bundle on $M$ defined by \eqref{e:sections} with Hermitian connection defined by \eqref{def:A}. In this case, the set  $\{\nu_{N,j}, j=0,1,2,\ldots\}$ of the eigenvalues of $\Delta^{L^N}$ is the union of two parts $\{\nu^{(i)}_{N,k}, k=0,1,\ldots,N-1\}$ and $\{\nu^{(c)}_{N,\ell}, \ell=0,1,2,\ldots\}$ given by \eqref{e:nu-i} and \eqref{e:nu-c}, respectively  (taking into account the multiplicities).
Denote by $\hat\varphi$ the Fourier transform of $\varphi$. 

\begin{theorem}[\cite{KT1}, Theorem 7]
For any $\varphi\in \mathcal S(\RR)$ and $1<E<\sqrt{2}$, i.e., $0 < E_0 = |p|^2 < 1$, one has an asymptotic expansion
\[
Y_N(\varphi) \sim \sum_{j=0}^\infty c_j(N,\varphi)N^{1-j}, \quad N\to \infty,
\]
where the coefficients $c_j(N,\varphi)$ are bounded in $N$. 

The coefficients $c_j$ can be computed explicitly. For the first two of them, we get
\begin{equation*}
\begin{aligned}
c_0(N,\varphi) = & (2g-2)E \hat\varphi (0)\\ & + (2g-2)E \sum_{k\neq 0} \hat\varphi \left(\frac{2\pi kE}{\sqrt{2-E^2}}\right) \exp(ik\pi)\exp\left(2\pi ik \sqrt{2-E^2}N\right), 
\end{aligned}
\end{equation*}
\begin{equation*}
\begin{aligned}
c_1(N,\varphi) =& (2g-2)2i \hat\varphi^\prime (0)  + \Bigg[\sum_{k\neq 0}(2g-2)2i \hat\varphi^\prime \left(\frac{2\pi kE}{\sqrt{2-E^2}}\right) \\
&+\sum_{k\neq 0} (2g-2)\frac{\pi ikE}{4\sqrt{2-E^2}}\hat\varphi \left(\frac{2\pi kE}{\sqrt{2-E^2}}\right)\\
& +\sum_{k\neq 0} (2g-2)i\frac{2\pi ikE}{(2-E^2)^{3/2}}\hat\varphi^{\prime\prime} \left(\frac{2\pi kE}{\sqrt{2-E^2}}\right)\Bigg] \\
& \ \ \ \ \times \exp(ik\pi)\exp\left(2\pi ik \sqrt{2-E^2}N\right).
\end{aligned}
\end{equation*}
\end{theorem}

In this paper we complete the computation of the Guillemin-Uribe trace formula for this example, considering the case $E\geq \sqrt{2}$. 

\begin{thm}
For any $\varphi\in \mathcal S(\RR)$ with compactly supported Fourier transform and $E\geq \sqrt{2}$, i.e., $E_0 = |p|^2 > 1$, one has an asymptotic expansion
\begin{equation}\label{e:YN2}
Y_N(\varphi) \sim \sum_{j=0}^\infty c_j(N,\varphi)N^{1-j}, \quad N\to \infty,
\end{equation}
where the coefficients $c_j(N,\varphi)$ are bounded in $N$. 

We have
\begin{equation}\label{e:c0}
c_0(N,\varphi) = (2g-2)E \hat\varphi (0),
\end{equation}
and, if ${\rm supp}\,\hat\varphi\subset \RR\setminus \{0\}$, then for $E>\sqrt{2}$,
\begin{multline}\label{e:c1a}
c_1(N,\varphi)=\sum_{h\in \{\Gamma\}_p}\sum_{k\neq 0} \frac{\log N(h)}{2\pi (N(h)^{k/2}-N(h)^{-k/2})}\frac{E}{\sqrt{E^2-2}}  \\
\times  \hat \varphi\left(\frac{E}{\sqrt{E^2-2}}k\log N(h) \right)\exp\left(-ik\log N(h) \sqrt{E^2-2}N\right),
\end{multline}
where $\{\Gamma\}_p$ is the set of representative of primitive conjugacy classes in $\Gamma$ and $N(h)$ stands for the norm of $h$ (see below for the definition),  
and for $E=\sqrt{2}$, 
\begin{equation}\label{e:c1b}
c_j(N,\varphi)=0, \quad j=1,2,\ldots. 
\end{equation}
\end{thm}

We give two proofs of this theorem. The first proof use directly the general Guillemin-Uribe formula and the description of the magnetic geodesic flow given in Section~\ref{s:Lie}. In the second proof, we use the results of Section~\ref{s:spectrum} to reduce our considerations in the case $E>\sqrt{2}$ to a spectral problem for the scaled Laplace-Beltrami operator, where we apply a version of the Guillemin-Uribe trace formula for the Laplace-Beltrami operator. 

\subsection{The case $E\geq\sqrt{2}$: Direct computation} 

Since all periodic trajectories of the magnetic geodesic flow  $\phi$ are non-degenerate, the existence and the form of the asymptotic expansion \eqref{e:YN2} follow from the general Guillemin-Uribe formula. It remains to compute the coefficients. Each coefficient is represented as an infinite sum, and each term of the sum corresponds either to $0$ or to a periodic trajectory. The contribution of $0$ to $c_0(N,\varphi)$ is given by
\begin{equation} \label{e:c0-gen0}
c^{(0)}_0(N,\varphi)=(2\pi)^{-2}\hat{\varphi}(0){\rm Vol}(X_E).
\end{equation}
In the current setting, the same computation as in \cite{KT1} gives \eqref{e:c0}. 

Let $E>\sqrt{2}$. Since each periodic trajectory $\gamma\subset X_E$ is non-degenerate (i.e. the map $I-P_\gamma$ is invertible, where $P_\gamma$ denotes the linear Poincar\'e map of $\gamma$), the contributions $c^{(\gamma)}_0(N,\varphi)$ and $c^{(\gamma)}_1(N,\varphi)$ of $\gamma$ into the formulas for the leading coefficient $c_0(N,\varphi)$ and the next term $c_1(N,\varphi)$, respectively, are given by
\begin{equation} \label{e:c0-nondegenerate}
c^{(\gamma)}_0(N,\varphi)=0, \quad c^{(\gamma)}_1(N,\varphi)=\frac{T^\#_\gamma e^{\pi im_\gamma/4}}{2\pi |\det(I-P_\gamma)|^{1/2}}e^{-iNS_\gamma}\hat\varphi(T_\gamma)
\end{equation}
where $T^\#_\gamma$ and $m_\gamma$ are  the primitive period and the Maslov index of $\gamma$, respectively. Note that here we use a slightly different notation than in \cite{KT1}. 

Recall that the lift of the magnetic geodesic flow on $X_E$ to $T^*\mathbb H$ is isomorphic to the flow $\Phi^{\alpha,\beta}_t$ on $S\mathbb H$ given by the right translation by $\exp(t(\alpha E_1+\beta E_3))$ with $\alpha$ and $\beta$ given by \eqref{ab}. Therefore, the lift of a periodic trajectory $\gamma\subset X_E$ is given by the curve $\{g\exp(t(\alpha E_1+\beta E_3)), t\in [0,T]\}$ such that 
\begin{equation} \label{e:gamma}
g\exp(T(\alpha E_1+\beta E_3))=h g
\end{equation} 
with some $T$ and $h=\begin{pmatrix} a & b\\ c & d\end{pmatrix}\in \Gamma$. It is clear that $\gamma$ depends only on the conjugacy class $\{h\}_\Gamma$ of $h$ in $\Gamma$, and the trajectory $\gamma$ is primitive if and only if the conjugacy class $\{h\}_\Gamma$ is primitive. One can show that for any primitive conjugacy class $\{h\}_\Gamma$ in $\Gamma$, there exists a unique primitive periodic trajectory $\gamma\subset X_E$, which satisfies \eqref{e:gamma} with some $h\in \{h\}_\Gamma$. 

Since the group $\Gamma$ is cocompact, each element $h\in \Gamma\setminus\{e\}$ is hyperbolic, $\operatorname{tr} h>2$. Therefore, it is conjugate to a unique element of the form  $\begin{pmatrix} N(h)^{1/2} & 0\\ 0 & N(h)^{-1/2}\end{pmatrix}$ with some $N(h)>1$, called the norm of $h$. 

On the other hand, we know that the flow $\Phi^{\alpha,\beta}_t$ is conjugate to the flow $\Phi^{\delta,0}_t$ with $\delta=\frac{\sqrt{E^2 - 2}}{E}$. Therefore, we get
\[
\exp(T^\#_\gamma\delta E_1)=\begin{pmatrix} e^{T^\#_\gamma\delta/2} & 0\\ 0 & e^{-T^\#_\gamma\delta/2}\end{pmatrix} =\begin{pmatrix} N(h)^{1/2} & 0\\ 0 & N(h)^{-1/2}\end{pmatrix},
\]
which gives $T^\#_\gamma=\frac{E}{\sqrt{E^2 - 2}}\log N(h)$.

We can write the period $T$ as $T=kT^\#_\gamma$ with some $k\in \ZZ\setminus \{0\}$. The Poincar\'e map $P_\gamma$ of $\gamma$ is a linear hyperbolic map with the eigenvalues $e^{\delta T}=N(h)^k$ and $e^{-\delta T}=N(h)^{-k}$. Therefore, we have 
\[
|\det(I-P_\gamma)|^{1/2}=((1-N(h)^k)(1-N(h)^{-k}))^{1/2}=|N(h)^{k/2}-N(h)^{-k/2}|.
\]
This also implies that 
\[
m_\gamma=0.
\]

Denote by $h_A(\gamma)\in S^1=\RR/2\pi \ZZ$ the holonomy of the projection $\pi_M\circ \gamma$ of the curve $\gamma$ to $M$ with respect to the connection $\nabla^L$ on $L$. Then the action $S_\gamma$ of $\gamma$ is defined modulo multiplies of $2\pi$ and given by (see \cite{KT1})
\[
S_\gamma=\frac{E^2-1}{E}T+ h_A(\gamma).
\]

To compute the action $S_\gamma$ of the periodic trajectory $\gamma$, we lift $\gamma$ to $S\mathbb H$ as above. The curve $\gamma$ on $S\mathbb H$ is not closed and we should use the formula \eqref{e:sections} to get a correct identification of the fibers of the line bundle $L$ at its extreme points. Since the form $F$ is exact on $\mathbb H$, $F=dA$ with $A$ given by \eqref{def:A}, by \eqref{e:gamma} and \eqref{e:sections}, we have 
\[
h_A(\gamma)=\int_{\pi_M\circ \gamma} A-2\arg(cz+d).
\]
Using \eqref{SH}, we compute
\[
\int_{\pi_M\circ \gamma} A=\int_0^T\frac{\dot x}{y}dt=-\int_0^T\alpha \sin\varphi dt=-\int_0^T(\dot \varphi-\beta)dt=\varphi(0)-\varphi(T)-\frac{1}{E}T. 
\]
Using \eqref{e:gamma} and computing the action of $h$ on $S\mathbb H$, we get with $z=x+iy$ 
\[
z(T)=\frac{az(0)+b}{cz(0)+d}, \quad \varphi(T)=\varphi(0)-2\arg(cz+d).
\]
We conclude that 
\[
h_A(\gamma)=-\frac{1}{E}T
\]
and
\[
S_\gamma= \frac{E^2-1}{E}kT^\#_\gamma-\frac{1}{E}kT^\#_\gamma=k\log N(h) \sqrt{E^2-2},
\]
This completes the proof of \eqref{e:c1a}.

If $E=\sqrt{2}$, then $E_0 = 1$ and the magnetic geodesic flow on $X_E$ is isomorphic to the horocyclic flow. It is well-known that this flow has no periodic trajectories, which implies \eqref{e:c1b}. 

\subsection{The case $E>\sqrt{2}$: Reduction to the Laplace-Beltrami operator} 
In this section, we use the relation \eqref{e:nu-c} to reduce our considerations in the case $E>\sqrt{2}$ to a spectral problem for the scaled Laplace-Beltrami operator $\Delta_M^{(0)}$, where we apply a version of a semiclassical trace formula.  

First, we write 
\[
Y_N(\varphi)=Y^{(i)}_N(\varphi)+Y^{(c)}_N(\varphi), 
\]
where 
\[
Y^{(i)}_N(\varphi)=\sum_{k=0}^{N-1}m_{N,k}\varphi \left(\sqrt{\nu_{N,k}^{(i)}+N^2}-EN\right),
\]
and 
\[
Y^{(c)}_N(\varphi)=\sum_{\ell=0}^{\infty}\varphi \left(\sqrt{\nu_{N,\ell}^{(c)}+N^2}-EN\right).
\]
Since $\nu_{N,k}^{(i)}\leq N^2$ for any $N\in \NN$ and $k=0,\ldots,N-1$, it is easy to see that $Y^{(i)}_N(\varphi)=O(N^{-\infty})$ as $N\to \infty$.

Using \eqref{e:nu-c}, we get
\[
Y^{(c)}_N(\varphi)=\sum_{\ell=0}^{\infty}\varphi\left(\sqrt{\lambda_\ell+2N^2}-EN\right).
\]

The right hand side of the last formula is closely related to the semiclassical trace formula for the Schr\"odinger operator. We compute its asymptotic expansion by applying Guillemin-Uribe trace formula in the following setting (cf. \cite[Corollary 7.5]{Gu-Uribe89}). Let $\Delta_M^{(0)}$ be the Laplace--Beltrami operator on $M$ associated with the Riemannian metric 
\[
g^{(0)}= \frac{2}{y^2}(dx^2+dy^2).
\]   
It is clear that $\Delta_M^{(0)}=\frac 12\Delta_M$. We will interpret this operator as the magnetic Laplacian associated with the vanishing magnetic field $F^{(0)}=0$. So the associated Hermitian line bundle $L_0$ is trivial, the Hermitian connection $\nabla^{L_0}$ is trivial and the associated magnetic Laplacian $\Delta^{L_0^N}$ coincides with $\Delta_M^{(0)}$ for all $N$. The eigenvalues $\nu^{(0)}_{N,\ell}$ of $\Delta^{L_0^N}$ are given by 
\[
\nu^{(0)}_{N,\ell}=\frac{1}{2}\lambda_\ell, \quad \ell=0,1,2,\ldots.
\]
For an arbitrary function $\psi\in \mathcal S(\RR)$, the sequence $Y^{(0)}_N(\psi)$ associated with the operator $\Delta^{L_0^N}$ and an energy level $E^{(0)}>1$ by the formula \eqref{e:Yp} has the form
\[
Y^{(0)}_N(\psi)=\sum_{\ell=0}^{\infty}\psi\left(\sqrt{\nu^{(0)}_{N,\ell}+N^2}-E^{(0)}N\right),\quad N\in \mathbb N.
\]
It is easy to see that 
\[
Y^{(c)}_N(\varphi)=Y^{(0)}_N(\psi)
\]
with 
\[
E^{(0)}=\frac{1}{\sqrt{2}}E,\quad \psi(z)=\varphi\left(\sqrt{2} z\right).
\]

By \cite[Corollary 7.5]{Gu-Uribe89}, the sequence $Y^{(0)}_N(\psi)$ admits an asymptotic expansion
\begin{equation*}
Y^{(0)}_N(\psi)\sim \sum_{j=0}^\infty c^{(0)}_j(N,\psi)N^{1-j},\quad N\to \infty,
\end{equation*}
where the coefficients $c^{(0)}_j(N,\psi)$ are bounded in $N$.  

The associated magnetic geodesic flow is the geodesic flow of $g^{(0)}$, that is, the Hamiltonian flow defined by the Hamiltonian 
\begin{equation*}
H^{(0)}(x,y,p_x,p_y)=\left(\frac{y^2}{2}(p_x^2+p_y^2)+1\right)^{1/2}
\end{equation*}
on the cotangent bundle $X=T^*M$ equipped with the standard symplectic form.
Put 
\[
X^{(0)}_{E^{(0)}}=(H^{(0)})^{-1}(E^{(0)}).
\]

As in \eqref{e:c0-gen0}, the contribution of $0$ is given by
\[
c^{(0)}_0(N,\psi)=(2\pi)^{-2}\hat{\psi}(0){\rm Vol}(X^{(0)}_{E^{(0)}}).
\]
We compute  
\[
{\rm Vol}_{g^{(0)}}(M)= 2\mathrm{Vol}(M)=4\pi(2g-2)
\]
and
\[
{\rm Vol}(X^{(0)}_{E^{(0)}})=2\pi E^{(0)}{\rm Vol}_{g^{(0)}}(M)=(2\pi)^2 (2g-2)E\sqrt{2}. 
\]
Finally, we observe that
\[
\hat \psi(z)=\frac{1}{\sqrt{2}} \hat \varphi\left(\frac{z}{\sqrt{2}}\right),
\]
Taking all this into account, we conclude that 
\[
c_0(N,\varphi)=c^{(0)}_0(N,\psi)=(2g-2)E\hat \varphi(0).
\]

As in \eqref{e:c0-nondegenerate}, the contributions of a periodic trajectory $\gamma\subset X_{E^{(0)}}$ with period $T_\gamma=kT^\#_\gamma$ are given by
\[
c^{(0,\gamma)}_0(N,\psi)=0, \quad 
c^{(0,\gamma)}_1(N,\psi)=\frac{T^\#_\gamma e^{\pi im_\gamma/4}}{2\pi |\det(I-P_\gamma)|^{1/2}}e^{-iNS_\gamma}\hat\psi(T_\gamma),
\]
where $P_\gamma$ denotes the Poincar\'e map of $\gamma$, $T^\#_\gamma$ and $m_\gamma$ are  the primitive period and the Maslov index of $\gamma$, respectively.

The Hamiltonian flow of $H^{(0)}$ on $X$ is given by
\[
\dot x=\frac{y^2}{2H^{(0)}}p_x,\quad \dot y=\frac{y^2}{2H^{(0)}}p_y,\quad \dot p_x=0,\quad \dot p_y=-\frac{y}{2H^{(0)}}(p_x^2+p_y^2),
\]
and its restriction to $X_{E^{(0)}}$ by
\[
\dot x=\frac{y^2}{2E^{(0)}}p_x,\quad \dot y=\frac{y^2}{2E^{(0)}}p_y,\quad \dot p_x=0,\quad \dot p_y=-\frac{y}{2E^{(0)}}(p_x^2+p_y^2).
\]
We introduce on $X_{E^{(0)}}= \{y^2(p_x^2 + p_y^2) = 2 ((E^{(0)})^2-1) \}$
the coordinates
$(x,y,\theta)$: 
\[
p_x= \frac{\sqrt{2 ((E^{(0)})^2-1)}}{y}\cos\theta, \quad p_y= \frac{\sqrt{2 ((E^{(0)})^2-1)}}{y}\sin\theta
\]
in which the system takes the form
\[
\dot x=\delta_0 y\cos\theta,
\quad \dot y=\delta_0 y\sin\theta,
\quad \dot \theta=-\delta_0\cos\theta,
\]
with 
\[
\delta_0=\frac{\sqrt{2 ((E^{(0)})^2-1)}}{2E^{(0)}}
\] 
and defines the flow $\Phi^{\delta_0,0}_t$.

As above, we get
\[
T^\#_\gamma=\frac{1}{\delta_0}\log N(h)=\frac{2E^{(0)}}{\sqrt{2 ((E^{(0)})^2-1)}}\log N(h)=\frac{E\sqrt{2}}{\sqrt{E^2-2}}\log N(h). 
\]
\[
|\det(I-P_\gamma)|^{1/2}=|N(h)^{k/2}-N(h)^{-k/2}|, \quad 
m_\gamma=0.
\]
The action $S_\gamma$ of $\gamma$ is given by
\[
S_\gamma=\frac{(E^{(0)})^2-1}{E^{(0)}}kT^\#_\gamma=k \log N(h) \sqrt{E^2-2}.
\]

We arrive at the desired formula:
\begin{multline*}
c^{(\gamma)}_1(N,\varphi)=c^{(0,\gamma)}_1(N,\psi)\\ =\frac{\log N(h)}{2\pi |N(h)^{k/2}-N(h)^{-k/2}|}\frac{E}{\sqrt{E^2-2}}  \hat \varphi\left(\frac{E}{\sqrt{E^2-2}}k\log N(h) \right) \\
\times \exp\left(-ik\log N(h)\sqrt{E^2-2}N\right). 
\end{multline*}

\end{document}